\documentclass[12pt,reqno]{amsart}
\usepackage{enumerate, latexsym, amsmath, amsfonts, amssymb, amsthm, color}
\def\pmod #1{\ ({\rm{mod}}\ #1)}
\def\Z{\Bbb Z}
\def\N{\Bbb N}

\def\Q{\Bbb Q}

\def\l{\left}
\def\r{\right}
\def\bg{\bigg}
\def\({\bg(}
\def\){\bg)}
\def\t{\text}
\def\f{\frac}

\def\ls{\leqslant}
\def\gs{\geqslant}
\def\se {\subseteq}

\def\bi{\binom}

\def\ve{\varepsilon}

\def\eq{\equiv}

\def\Proof{\noindent{\it Proof}}

\theoremstyle{plain}
\newtheorem{theorem}{Theorem}

\newtheorem{lemma}{Lemma}
\newtheorem{corollary}{Corollary}

\theoremstyle{definition}

\theoremstyle{remark}
\newtheorem{remark}{Remark}

 \vspace{4mm}

\begin{document}

\hbox{Accepted for publication in Proc. of the 2019 Asian Logic Conf. (World Sci.)}
\medskip

\title
[{On Diophantine equations over $\mathbb Z[i]$ with $52$ unknowns}]
{On Diophantine equations over $\mathbb Z[i]$\\ with $52$ unknowns}

\author
[Yu. Matiyasevich and Z.-W. Sun] {Yuri Matiyasevich and Zhi-Wei Sun}

\address {(Yuri Matiyasevich) St. Petersburg Department of Steklov Mathematical Institute of Russian Academy of Sciences, Fontanka 27, 191023, St. Petersburg, Russia}
\email{yumat@pdmi.ras.ru}

\address{(Zhi-Wei Sun, corresponding author) Department of Mathematics, Nanjing
University, Nanjing 210093, People's Republic of China}
\email{zwsun@nju.edu.cn}

\keywords{Hilber's Tenth Problem, Diophantine equation, Gaussian ring, undecidability.
\newline \indent The work is supported by the NSFC-RFBR Cooperation and Exchange Program (grants NSFC 11811530072 and
RFBR 18-51-53020-GFEN-a). The second author is also supported by the Natural Science Foundation of China (grant no. 11971222).}

\begin{abstract} In this paper we show that there is no algorithm to decide whether
an arbitrarily given polynomial equation $P(z_1,\ldots,z_{52})$ $=0$
(with integer coefficients) over the Gaussian ring $\mathbb Z[i]$
is solvable.
\end{abstract}
\maketitle

\section{Introduction}
\setcounter{lemma}{0}
\setcounter{theorem}{0}
\setcounter{corollary}{0}
\setcounter{remark}{0}

The original HTP (Hilbert's Tenth Problem) asks for an (effective) algorithm  to test whether an arbitrary polynomial Diophantine equation
with integer coefficients has solutions over the ring $\Z$ of the integers. This was finally solved by Yu. Matiyasevich \cite{M70} negatively in 1970 based on the work of M. Davis, H. Putnam and J. Robinson
\cite{DPR}. Z.-W. Sun \cite{S17} showed further that there is no algorithm
to decide for any given $P(x_1,\ldots,x_{11})\in\Z[x_1,\ldots,x_{11}]$ whether the equation
$P(x_1,\ldots,x_{11})=0$ has integer solutions.

Let $K$ be a number field which is a finite extension of the field $\Q$ of rational numbers.
It is natural to ask whether HTP over the ring $O_K$ of algebraic integers in $K$
is unsolvable. Clearly, if $\Z$ is Diophantine over $O_K$ then
HTP over $O_K$ is undecidable with the aid of Matiyasevich's theorem.
It is known that $\Z$ is Diophantine over $O_K$ if $[K:\Q]=2$ or $K$ is totally real (cf. J. Denef \cite{D75,D80}),
or $[K:\Q]\gs3$ and $K$ has exactly two nonreal embeddings into the field of complex numbers
(cf. T. Pheidas \cite{P88}), or $K$ is an abelian number field (cf. H. N. Shapiro and A. Shlapentokh \cite{SS}).

In this paper we study Diophantine equations with few unknowns over the Gaussian ring
$$\Z[i]=O_{\Q(i)}=\{a+bi:\ a,b\in\Z\}.$$
Our main results are as follows.

\begin{theorem}\label{Th1.1} A number $z\in\Z[i]$ is a rational integer if and only if
there are $v,w,x,y\in\Z[i]$ with $v\not=0$ such that
\begin{equation}\label{zzi}
\begin{aligned}&(4(2v(2(2z+1)^2+1)-y)^2-3y^2-1)^2
\\&+2(w^2-1-3y^2(2z+1-xy)^2)^2=0.\end{aligned}
\end{equation}
\end{theorem}

\begin{theorem}\label{Th1.2} For any r.e. (recursively enumerable) set $\mathcal A\se\N=\{0,1,2,\ldots\}$, there is a polynomial
$P(z_0,z_1,\ldots,z_{52})$ with integer coefficients such that
for any $a\in\N$ we have
\begin{equation}\label{re}a\in \mathcal A\iff P(a,z_1,\ldots,z_{52})=0\ \t{for some}\ z_1,\ldots,z_{52}\in\Z[i].\end{equation}
\end{theorem}

It is well known (cf. N. Cutland \cite{C80}) that there are nonrecursive r.e. subsets of $\N$. Thus Theorem \ref{Th1.2}
has the following corollary.

\begin{corollary}\label{Cor1.1} There is no algorithm to decide for any polynomial
$P(z_1,\ldots,z_{52})$ with integer coefficients whether the equation
$$P(z_1,\ldots,z_{52})=0$$
has solutions in $\Z[i]$.
\end{corollary}

We will provide some lemmas in the next section and then show Theorems 1.1-1.2 in Section 3.

\section{Some Lemmas}
\setcounter{lemma}{0}
\setcounter{theorem}{0}
\setcounter{corollary}{0}
\setcounter{remark}{0}

For $A,B\in\Z$, the Lucas sequence $(u_n(A,B))_{n\gs0}$ is given by
$u_0(A,B)=0$, $u_1(A,B)=1$,
and $$\ u_{n+1}(A,B)=Au_n(A,B)-Bu_{n-1}(A,B)\ \ (n=1,2,3,\ldots).$$
Sun \cite{S92} studied arithmetic properties of such sequences as well as related Diophantine representations over $\Z$.

\begin{lemma}\label{Lem2.1} Let $A,B\in\Z$.

{\rm (i)} For any $k,n,r\in\N$, we have the identity
$$u_{kn+r}(A,B)=\sum_{j=0}^n\bi nj(u_{k+1}(A,B)-Au_k(A,B))^{n-j}u_k^ju_{j+r}.$$

{\rm (ii)} Let $A,B,M\in\Z$ with $M\not=0$. Then $B$ is relatively prime to $M$
if and only if
$u_n(A,B)\eq0\pmod M$ and $u_{n+1}(A,B)\eq1\pmod M$ for some $n\in\Z^+=\{1,2,3,\ldots\}$.

{\rm (iii)} If $A>B\gs0$, then $(A-B)^n\ls u_{n+1}(A,B)\ls A^n$ for all $n\in\N$.
\end{lemma}
\begin{remark} Parts (i)-(iii) are Lemmas 2, 6, 8 of Sun \cite{S92}.
\end{remark}

\begin{lemma}\label{Lem2.2} Let $A\in\{2,3,\ldots\}$. Then
$$x^2-Axy+y^2=1\ \t{with}\ x,y\in\N\ \t{and}\ x\gs y$$
if and only if
$$x=u_{n+1}(A,1)\ \t{and}\ y=u_n(A,1)\ \t{for some}\ n\in\N.$$
\end{lemma}
\begin{remark} This is a known result, see, e.g., Sun \cite[Lemma 9]{S92}.
\end{remark}
\begin{lemma}\label{Lem2.3} If $x,y\in\Z[i]$ and
$x^2-4xy+y^2=1$, then $x,y\in\Z$.
\end{lemma}
\begin{remark}  This follows from a more general result of Denef
\cite{D75}; a proof for this particular case was also presented in
Matiyasevich \cite[Section 7.3]{M93}.
\end{remark}

\begin{lemma}\label{Lem2.4} For $x,y\in\Z[i]$, we have
$$x=0\land y=0\iff x^2+2y^2=0.$$
\end{lemma}
\Proof. Though the result is known, here we provide a simple proof.

 Suppose that $x^2+2y^2=0$ but $x\not=0$ or $y\not=0$. Then $xy\not=0$ and $x/y\in\{\sqrt2\,i,-\sqrt2\,i\}$. As $x/y\in\Q(i)=\{r+si:\ r,s\in\Q\}$, and $\sqrt2$
is irrational, we obtain a contradiction. This ends the proof. \qed

\begin{lemma}\label{Lem2.5} An integer $m$ is nonzero if and only if $m=(2x+1)(3y+1)$
for some $x,y\in\Z$.
\end{lemma}
\begin{remark} This is a useful observation of S.-P. Tung \cite{T85}.
\end{remark}

\section{Proofs of Theorems 1.1 and 1.2}
\setcounter{lemma}{0}
\setcounter{theorem}{0}
\setcounter{corollary}{0}
\setcounter{remark}{0}

\medskip
\noindent{\it Proof of Theorem 1.1}. (i) We first show the ``if" direction.

Suppose that there are $v,w,x,y\in\Z[i]$ with $v\not=0$ satisfying \eqref{zzi}.
In view of Lemma \ref{Lem2.4}, we have
\begin{equation}\label{2.1}4(2v(2(2z+1)^2+1)-y)^2-3y^2-1=0
\end{equation}
and
\begin{equation}\label{2.2}
w^2-1-3y^2(2z+1-xy)^2=0.
\end{equation}
Let $y_*=4v(2(2z+1)^2+1)$ and $w_*=w+2(2z+1-xy)y$. Then
$$y_*^2-4y_*y+y^2=(y_*-2y)^2-3y^2=1$$
and
\begin{align*}&w_*^2-4w_*y(2z+1-xy)+y^2(2z+1-xy)^2
\\=&(w_*-2y(2z+1-xy))^2-3y^2(2z+1-xy)^2
\\=&w^2-3y^2(2z+1-xy)^2=1.\end{align*}
Applying Lemma \ref{Lem2.3}, we see that $y,y_*,w_*,y(2z+1-xy)\in\Z$.
Thus both $2z+1-xy$ and $w$ are rational integers.

Note that
$$\f{|y_*|}4\gs2|2z+1|^2-1=|2z+1|(2|2z+1|-1)+(|2z+1|-1)\gs|2z+1|$$
and
$$(y-2y_*)^2=3y_*^2+1\ls 3y_*^2+\f{y_*^2}{16}=\l(\f74y_*\r)^2.$$
If $(y-2y_*)^2=(\f 74y_*)^2$, then we must have $|y_*|/4=1=|2z+1|$, hence $z\in\{0,-1\}$ and $|y_*|=|12v|>4$.
Therefore
\begin{align*}|y|>2|y_*|-\f 74|y_*|=\f{|y_*|}4\gs|2z+1|.
\end{align*}
Recall that $2z+1-xy\in\Z$, and write $x=a+bi$ with $a,b\in\Z$. Then
$$|y|^2>|2z+1|^2=|(2z+1-xy)+(a+bi)y|^2=(2z+1-xy+ay)^2+b^2y^2,$$
hence $b=0$ and $x\in\Z$. Thus $2z+1\in\Z$ and hence $z\in\Z$.
\medskip

(ii) Below we show the ``only if" direction. For $n\in\N$ we simply write $u_n$ to denote $u_n(4,1)$.

Let $z\in\Z$ and $k=|2z+1|$. By Lemma \ref{Lem2.1}(ii), for some $n\in\N$ we have
$u_{n+1}\eq0\pmod{4(2k^2+1)}$. In view of Lemma \ref{Lem2.1}(iii), $u_{kn}\gs 3^{kn-1}$
and $u_{n+1}\gs 3^{n}$.
Write $u_{n+1}=4(2k^2+1)v$ with $v\in\Z^+$
and set $y=u_n$. Then
$$4(2v(2k^2+1)-y)^2=(u_{n+1}-2u_n)^2=3u_n^2+1=3y^2+1$$
with the aid of Lemma \ref{Lem2.2}.
By Lemma \ref{Lem2.1}(i),
$$u_{nk}\eq k(u_{n+1}-4u_n)^{k-1}u_n\pmod{u_n^2}.$$
Let $q=u_{kn}/u_n\in\Z^+$. Then
$$q\eq ku_{n+1}^{k-1}\eq k\pmod {u_n}$$
since $k\eq1\pmod2$ and $u_{n+1}^2=1-u_n^2+4u_nu_{n+1}\eq 1\pmod{u_n}$.
Define $\ve=1$ if $z\gs0$, and $\ve=-1$ if $z<0$. Then
$\ve u_{kn}=u_n(\ve k-xu_n)=y(2z+1-xy)$ for some $x\in\Z$.
Let $w_*=\ve u_{kn+1}$ and $w=w_*-2\ve u_{kn}$. Then
$$w^2-3y^2(2z+1-xy)^2=(u_{kn+1}-2u_{kn})^2-3u_{kn}^2=1$$
by Lemma \ref{Lem2.2}.
Now it is clear that \eqref{zzi} holds.

In view of the above, we have completed the proof of Theorem \ref{Th1.1}. \qed

\begin{remark} In view of Lemma \ref{Lem2.5} and the proof of Theorem 1.1,
a number $z\in\Z[i]$ is an rational integer if and only if there are $s,t,w,x,y\in\Z[i]$
such that \eqref{zzi} holds with $v=(2s+1)(3t+1)$.
\end{remark}

\medskip
\noindent{\it Proof of Theorem \ref{Th1.2}}. By Sun \cite[Theorem 1.1(ii)]{S17}, there is a
polynomial $f(z_0,\ldots,z_{10})\in\Z[z_0,\ldots,z_{10}]$ such that $a\in\N$ belongs to ${\mathcal A}$
if and only if $f(a,z_1,\ldots,z_{10})=0$ for some $z_1,\ldots,z_{10}\in\Z$ with $z_{10}\not=0$.

Let $F(v,w,x,y,z)$ denote the left-hand side of \eqref{zzi}.
For $z_k\in\Z[i]$, by Theorem \ref{Th1.1}, $z_k\in\Z$ if and only if
$F(v_k,w_k,x_k,y_k,z_k)=0$ for some $v_k,w_k,x_k,y_k\in\Z[i]$ with $v_k\not=0$.
Thus, $a\in\mathcal A$ if and only if there are
$$v_k,w_k,x_k,y_k,z_k\in\Z[i]\ \ (k=1,\ldots,10)$$
with $F(v_k,w_k,x_k,y_k,z_k)=0\ \t{for all}\ k=1,\ldots,10$
such that $z_{10}\prod_{k=1}^{10}v_k\not=0$ and $f(a,z_1,\ldots,z_{10})=0$.
By the proof of Theorem \ref{Th1.1}, when $a\in\mathcal A$ we can actually choose
$z_{10},v_1,\ldots,v_{10}\in\Z\setminus\{0\}$ to meet the requirements.
Therefore, in view of Lemma \ref{Lem2.5},
$a\in\mathcal A$ if and only if there are
$$v_k,w_k,x_k,y_k,z_k\in\Z[i]\ \ (k=1,\ldots,10)$$
such that $f(a,z_1,\ldots,z_{10})=0$,  $F(v_k,w_k,x_k,y_k,z_k)=0\ \t{for all}\ k=1,\ldots,10$,
and $z_{10}\prod_{k=1}^{10}v_k=(2s+1)(3t+1)$ for some $s,t\in\Z[i]$.
Thus, in light of Lemma \ref{Lem2.4}, \eqref{re} holds for some polynomial
$P(z_0,z_1,\ldots,z_{52})\in\Z[z_0,z_1,\ldots,z_{52}]$.
This concludes the proof. \qed

\end{document}